\documentclass[12pt]{article}
\usepackage{amsmath, amssymb, amsthm, enumerate} 
\usepackage[utf8]{inputenc}
\usepackage{mathtools}

\newtheorem{theorem}{Theorem}[section]

\newtheorem{lemma}[theorem]{Lemma}

\def\rr{{\mathbb R}}
\def\zz{{\mathbb Z}}

\def\nn{{\mathbb N}} 
\def\ff{{\cal F}}
\def\rn{\rr ^n}
\def\zn{\zz ^n}
 
\def\am{^{-1}} 
\def\su{\subset}

\def\al{\alpha}

\def\ga{\gamma}
\def\de{\delta}
\def\De{\Delta}
\def\ep{\varepsilon}

\def\la{\lambda}
\def\Om{\Omega}
\def\cd{\cdot}
\def\stb{,\ldots ,}

\def\dim{{\rm dim}\, }

\def\ol{\overline}
\def\V{\Vert}
\def\ran{\rangle}
\def\lan{\langle}
\def\am{^{-1}}
\def\nl{[0,1]}

\def\bsk{\bigskip}
\def\noi{\noindent}

\title{A superposition theorem of Kolmogorov type for bounded continuous
  functions}

\author{Mikl\'os Laczkovich\thanks{Partially supported by Hungarian Scientific
Foundation grant no. K124749.}}

\begin{document}

\maketitle 

\begin{abstract}
Let $C(\rn )$ denote the set of real valued continuous functions defined on
$\rn$. We prove that for every $n\ge 2$ there are
positive numbers $\la _1 \stb \la _n$ and continuous functions $\phi_1 \stb
\phi _m \in C(\rr )$ with the following property: for every bounded
and continuous $f\in C(\rn )$ there is a continuous function $g\in C(\rr )$
such that
$$f(x)=\sum_{q=1}^m g\left( \sum_{p=1}^n \la _p \phi _q (x_p ) \right)$$
for every $x=(x_1 \stb x_n )\in \rn$.
Consequently, every $f\in C(\rn )$ can be obtained from
continuous functions of one variable using compositions and additions.
\end{abstract}

\insert\footins{\footnotesize{MSC code: 26B40}}
\insert\footins{\footnotesize{Key words: Kolmogorov superposition theorem}}

\section{Introduction and main results}
By a well-known theorem of Kolmogorov, every continuous function defined on
$\nl ^n$ can be represented in the form $\sum_{q=1}^{2n+1}
h_q (\phi _{q,1} (x_1 )+\ldots +\phi _{q,n} (x_n ) )$, where $h_q$ and $\phi _{q,p}$
are continuous functions of one variable for every $q=1\stb 2n+1$ and $p=1
\stb n$ \cite{Ko}. The theorem was subsequently improved and its proof was
simplified; see \cite{M}, \cite{S} and \cite{V} for these developments.
One of the sharpest variants with a particularly simple proof was presented by
Kahane \cite{Ka}.

Kolmogorov's theorem was put into topological context by a series of papers.
The following concept was introduced along this line. Let $C(X)$ denote the
set of real valued continuous functions defined on a topological space $X$.
A family $\Phi \su C(X)$ is said to be basic (resp. basic*) if every continuous
function $f\in C(X)$ (resp. every bounded and continuous $f\in C(X)$)
can be represented as $f=\sum_{i=1}^m g_i \circ \phi _i$ for some
$\phi _i \in \Phi$ and $g_i \in C(\rr )$ $(i=1\stb m)$. 

It was proved by S. Demko that in $\rn$ there exists a basic* family
consisting of $2n+1$ functions \cite{D}. Then Y. Hattori proved that 
if $X$ is a locally compact separable metric space with $\dim X\le n$,
then there exists a basic family on $X$
consisting of $2n+1$ functions (see \cite{H} and \cite{FG}).
In particular, there is such a basic family on $\rn$.

In this note we prove that there is a basic* family on $\rn$ consisting of
functions of the form $\phi_1 (x_1 )+\ldots +\phi _n (x_n )$. Our main result
is the following.
\begin{theorem}\label{t1}
Let $n\ge 2$ and $m>(2+\sqrt 2 )n$ be integers, and let $\la _1 \stb \la _n$ be distinct positive numbers. Then there are continuous
functions $\phi_1 \stb \phi _m \in C(\rr )$ with the following property:
for every bounded $f\in C(\rn )$ there is a continuous function $g\in C(\rr )$
such that
\begin{equation}\label{e10}
f(x)=\sum_{q=1}^m g\left( \sum_{p=1}^n \la _p \phi _q (x_p ) \right)
\end{equation}
for every $x=(x_1 \stb x_n )\in \rn$.
\end{theorem}

We note the following obvious consequence of Theorem \ref{t1}: {\it
every $f\in C(\rn )$ can be obtained from continuous functions of one
variable using compositions and additions.} Indeed, if $f\in C(\rr )$, then
$\ol f =\arctan f$ has a representation as in \eqref{e10}, and $f=\tan \ol f$.

The number of terms appearing in the theorem is probably not sharp, and a more
precise estimate of the norms of the approximating functions could reduce
the number $(2+\sqrt 2 )n$. It would be interesting to decide if 
such a basic* family containing $2n+1$ functions exists or not. 
Another natural question concerns the existence of basic families of this
form. The constructions of basic families in \cite{H} and elsewhere restrict
the range of the functions of the family in such a way that seems to exclude
functions of the form $\phi_1 (x_1 )+\ldots + \phi _n (x_n )$. Again, it would
be interesting to decide if such basic families exist or not.

In all variants of Kolmogorov's theorem the inner functions can be chosen to be
monotonic. Our proof of Theorem \ref{t1} does not produce monotonic functions.
A modified construction gives such inner functions, but then the number
of terms must be increased. The following can be proved.
\begin{theorem}\label{t2}
Let $n\ge 2$ and $m>(2+\sqrt 2 )(2n-1)$ be integers. There are positive
numbers $\la _{q,p}$, and there are continuous and increasing
functions $\phi_1 \stb \phi _m \in C(\rr )$ with the following property:
for every bounded $f\in C(\rn )$ there is a continuous function $g\in C(\rr )$
such that
\begin{equation*}
f(x)=\sum_{q=1}^m g\left( \sum_{p=1}^n \la _{q,p} \phi _q (x_p ) \right)
\end{equation*}
for every $x=(x_1 \stb x_n )\in \rn$.
\end{theorem}

We prove Theorem \ref{t1} in the next three sections, and then sketch the proof
of Theorem \ref{t2} in Section 5.  

\section{Preliminaries}
We put $\V f\V _A =\sup_{x\in A} |f(x)|$ for every $f\in C(\rn )$ and $A\su \rn$.

Let $\tau _n$ denote the topology on $C(\rn )$ of uniform convergence on
compact sets. That is, let $G\in \tau _n$ if for every $f\in G$
there is a compact set $A\su \rn$ and there is a $\de >0$ such that $g\in G$
whenever $g\in C(\rn )$ and $\V g-f\V _A <\de$. Let
$$d(f,g) =\sum_{k=1}^\infty \min \left( \frac{1}{2^k} , \V f-g\V _{[-k ,k ]^n}
\right)$$
for every $f,g\in C(\rn )$. It is well-known that $d$ is a metric on $C(\rn )$
and the metric space $(C(\rn ),d)$ is complete and separable. The topology
generated by this metric space is $\tau _n$, hence $(C(\rn ), \tau _n )$ is 
a Polish space.  The set
$$\Phi _k =\{ f\in C(\rr ) \colon |f(x)-|x| |<1 \ (x\in [-k,k] )\}$$
is open in $C(\rr )$ for every $k>0$. Indeed, if $f\in \Phi _k$, then
there is a $\de >0$ such that $|f(x)-|x||<1-\de$ for every $x\in [-k,k]$, and
thus $g\in \Phi _k$ whenever $g\in C(\rr )$ and $\V g-f\V _{[-k,k]}<\de$.

Let $\Phi$ denote the set of continuous functions $f\colon \rr \to
\rr$ such that $|f(x)-|x||<1$ for every $x\in \rr$. Since $\Phi =
\bigcap_{k=1}^\infty \Phi _k$, it follows that $\Phi$ is a $G_\de$ subset of
$C(\rr )$. Therefore, $\Phi$ equipped with the subspace topology is a Polish
space (see \cite[Theorem 3.11]{K}). Then $\Phi ^{m}$, as a product space, is
also Polish.

In the sequel we fix the distinct positive numbers $\la _1 \stb \la _p$ such
that $\sum_{p=1}^n  \la _p =1$. The assumption of this extra condition does not
affect the proof, since, if the functions $\phi _1 \stb \phi _m$ satisfy the
requirement of Theorem 1, then the same is true for the functions
$c\cd \phi _1 \stb c\cd \phi _m$ for every $c\ne 0$. We put $\la =
\min_{1\le p\le n} \la _p$ and $C=1/\la$.

Let $\nn$ denote the set of nonnegative integers.
We put $D_0 =1$ and $D_{t+1}=C(D_t+10)$ $(t\in \nn )$. The cube
$[-D_t ,D_t ]^n$ will be denoted by $Q_t$. We denote $\ep _0 =(m-n)\am$ and
$\ep _1 =n\ep _0$.

Let $f\in C(\rn )$, $N\in \nn$ and $\eta >0$ be given.
We denote by $\Om _{N, \eta} (f)$ the set of those $m$-tuples
$(\phi _1 \stb \phi _{m})\in \Phi ^{m}$ for which there exists a function
$h\in C(\rr )$ such that
\begin{equation}\label{e1}
h(x)=0\ \text{if}\ x\le -2 \ \text{or} \ x\ge D_N +2 ,
\end{equation}
\begin{equation}\label{e2}
\V h\V _{[-D_t ,D_t ]} <  \ep _0 \V f\V _{Q_{t+1}} +\eta
\end{equation}
for every $t\in \nn$, and
\begin{equation}\label{e3}
\left\V f(x_1 \stb x_n ) -\sum_{q=1}^{m} h\left( \sum_{p=1}^n  \la _p
\phi _q (x_p ) \right) \right\V _{Q_t} < \ep _1 \V f\V _{Q_{t+1}} +\eta
\end{equation}
for every $0\le t\le N$.  
\
\section{Two lemmas}
\begin{lemma}\label{l1}
The set $\Om _{N,\eta} (f)$ is relatively open and dense in $\Phi ^m$ for every
$f\in C(\rn )$, $N\in \nn$ and $\eta >0$.
\end{lemma}
\proof First we show that $\Om _{N, \eta} (f)$ is relatively open. Let
$(\psi _1 \stb \psi _{m})\in \Om _{N, \eta} (f)$ be given, and let $h$ satisfy
\eqref{e1}-\eqref{e3} (with $\psi _q$ in place of $\phi _q$). 

Let $\ga >0$ be such that the difference between the two sides of \eqref{e3} is
greater than $\ga$ for every $0\le t\le N$.

Since $h$ is uniformly continuous by \eqref{e1}, there is a $\de >0$ such
that $|h(x)-h(y)|<\ga /m$ for every $x,y\in \rr$, $|x-y|<\de$.
Suppose that $(\phi _1 \stb \phi_m )\in \Phi ^m$ and $\V \phi _q -\psi _q
\V_{[-D_N ,D_N]} <\de$ $(q=1\stb m)$. We prove $(\phi _1 \stb \phi _{m})\in
\Om _{N,\eta} (f)$.

We put $X_q (x)=\sum_{p=1}^n \la _p \phi _q (x_p )$ and
$Y_q (x)=\sum_{p=1}^n \la _p \psi _q (x_p )$ for every $q=1\stb m$ and
$x=(x_1 \stb x_n )\in \rn$. If $x\in Q_N$, then
$$|X_q (x)- Y_q (x)|\le \sum_{p=1}^n \la _p |\phi _q (x_p )-
\psi _q (x_p )| <\de ,$$
and thus $|h(X_q (x))- h(Y_q (x))|<\ga /m$ for every $q$. Then
the left hand side of \eqref{e3} is increased by less than $\ga$ if we
replace $\psi _q$ by $\phi _q$. Therefore,
\eqref{e3} holds true, proving that $\Om _{N, \eta} (f)$ is open.

\bsk
Now we show that $\Om _{N, \eta} (f)$ is dense in $\Phi ^{m}$. Let
$(\psi _1 \stb \psi _{m})\in \Phi ^{m}$ and $M,\xi >0$ be given. We have to show
that $(\phi _1 \stb \phi _{m})\in \Om _{N,\eta} (f)$ for a suitable
$(\phi _1 \stb \phi _{m})\in \Phi ^{m}$ satisfying
$\V \phi _q - \psi _q \V _{[-M,M]} <\xi$ $(q=1\stb m)$. We may assume 
$\xi <(6n)\am$, and that
$$|x| -1<\psi _q (x)-\xi <\psi _q (x)+\xi <|x|+1$$
for every $x\in [-D_N ,D_N ]$ and $q=1\stb m$.

Since $f$ is uniformly continuous on $Q_N$ and $\psi _1 \stb \psi _m$
are uniformly continuous on $[-D_N ,D_N ]$, there is a $\de >0$ such that
\begin{enumerate}[{\rm(i)}]
\item $|f(x)-f(y)|<\eta /m$ whenever $x,y\in Q_N$ and 
$|x-y|<mn \de$, and
\item $|\psi _q (x)-\psi _q (y)|<\xi$ whenever $x,y\in [-D_N ,D_N ]$,
$|x-y|<m \de$ and $q=1\stb m$.
\end{enumerate}
We may assume $\de <1/(mn)$. Let
$$I_q (j)= [q\de +mj\de ,(q +m-1)\de + mj\de ]$$
for every $q=1\stb m$ and $j\in \zz$. Note that (a) for every fixed $q$,
the closed intervals $I_q (j)$ are pairwise disjoint, and (b) every $x\in \rr$
belongs to at least $m-1$ of the intervals $I_q (j)$.

For every $q$, the closed intervals $I_q (j) \cap [-D_N ,D_N ]$ are empty for
all but a finite number of integers $j$.  Let $K_q (1) \stb K_q (r_q )$ be the
enumeration of the nonempty intervals $I_q (j) \cap [-D_N ,D_N ]$ $(j\in \zz )$.

Since $|I_q (j)|=(m-1)\de$ for every $q$ and $j$, it follows from the
choice of $\de$ that the oscillation of $\psi _q$ on each of the intervals
$K_q (i)$ is less than $\xi$ for every $q$ and $i$. This implies that
there are continuous functions $\phi _1 \stb \phi _m \in
\Phi$ such that each $\phi _q$ is constant on each interval $K_q (i)$,
and $|\phi _q (x) - \psi _q (x)| <\xi$ for every $x\in \rr$ and $q=1\stb m$. 

We put $X_q (x_1 \stb x_n )=\sum_{p=1}^n \la _p \phi _q (x_p )$
$(x_1 \stb x_n  \in \rr )$. Then $X_q$ is constant on each box
$$B _q (i)= K_q (i_1 )\times \ldots \times K_q (i_n ) \quad 
(1\le i_1 \stb i_n \le r_q )$$
for every $q=1\stb m$. The system  $\{ B_q (i) \colon 1\le i_1 \stb i_n \le
r_q \}$ consists of pairwise disjoint boxes for every $q=1\stb m$.
Let $J$ denote the set $\{ (q,i)\colon q=1\stb m, \ i\in \{1\stb r_q \} ^n \}$.
Changing slightly the values of $\phi _q$ on the
intervals $K_q (i)$ we may assume that the values $u_{q,i} =X_q (B_q (i))$
$((q,i)\in J)$ are distinct. (This is made possible by the assumption that
the numbers $\la _1 \stb \la _n$ are distinct.) The proof of the lemma will be
complete if we show that $(\phi _1 \stb \phi _m )\in \Om _{N, \eta} (f)$.

For every $(q,i)\in J$ we choose an element $a_{q,i} \in B_q (i)$, and define
$h(u_{q,i})= \ep _0 f(a_{q,i} )$. Then we extend $h$ to $\rr$ as follows. Let
$u_1 <u_2 <\ldots < u_s$ be the
enumeration of the numbers $u_{q,i}$ $((q,i)\in J)$, where $s=|J|$. We put
$u_0 =u_1 -1$ and $u_{s+1} =u_s +1$. We define $h(x)=0$ for every $x\le u_0$ and
$x\ge u_{s+1}$, and define $h$ on the interval $[u_\nu ,u_{\nu +1}]$ linearly for
every $\nu =0\stb s$. The function $h$ defined this way is piecewise linear,
hence continuous on $\rr$. We prove \eqref{e1}-\eqref{e3}.

\noi
{\it Proof of} \eqref{e1} and \eqref{e2}:
For every $1\le \nu \le s$ we have $u_\nu = u_{q,i} =X_q (a_{q,i} )$ for some
$(q,i)\in J$. If $a_{q,i}=(a_{q,i,1}\stb a_{q,i,n})$, then we have
$$u_\nu = X_q (a_{q,i} ) =\sum_{p=1}^n \la _p \phi _q (a_{q,i,p} ) \le
\sum_{p=1}^n \la _p  (|a_{q,i,p}|+1)\le D_N +1,$$
and
$$u_\nu = X_q (a_{q,i} ) =\sum_{p=1}^n \la _p \phi _q (a_{q,i,p} ) \ge
\sum_{p=1}^n \la _p (|a_{q,i,p}|-1)\ge -1.$$
Thus $u_1 \ge -1$ and $u_s \le D_N +1$, and \eqref{e1} follows.

If $x\in [1 ,D_N ]$, then $y=(x \stb x)\in Q_N$. If $y\in B_{q,i}$,
then $|x -a_{q,i,p}|<m\de <1$ and $a_{q,i,p} >0$ for every $p$. Thus
\begin{align*}
|u_{q,i} -x| &= \left| X_q (a_{q,i} ) -x\right|  = \left| \sum_{p=1}^n \la _p
\phi _q (a_{q,i,p} ) -x\right| \le
\sum_{p=1}^n \la _p |\phi _q (a_{q,i,p}) -x| \\
&\le \sum_{p=1}^n \la _p |\phi _q (a_{q,i,p}) -a_{q,i,p}|+ \sum_{p=1}^n \la _p
|a_{q,i,p} -x| \le 1+ 1=2.
\end{align*}
That is, for every $x\in [1,D_N ]$ there is a $\nu$ such that 
$|x-u_\nu |\le 2$. Then, for every $x\in [-1,D_N +1 ]$ there is a $\nu$ such
that $|x-u_\nu |\le 4$. This implies $u_\nu -u_{\nu -1} \le 8$ for every
$1\le \nu \le s+1$. 

Next we prove that, for every $x\in \rr$, 
\begin{equation}\label{e7}
|h(x)| \le \ep _0 \V f\V _{T_x} , \ \text{where} \ T_x =[-C|x|-9C,C|x|+9C]^n .
\end{equation}
(Note that $C=1/\la$, where $\la =\min_{1\le p\le n} \la _p$.)
This is clear if $x\le u_0$ or $x\ge u_{s+1}$, so we may
assume $u_{\nu -1} \le x\le u_\nu$ for some $1\le \nu \le s+1$.
Then, by $u_\nu -u_{\nu -1} \le 8$ we have $|u_{\nu -1}|, |u_\nu | \le |x|+8$.

If $u_\nu = u_{q,i} =X_q (a_{q,i} )$, then
\begin{align*}
|a_{q,i}| &\le \sum_{p=1}^n |a_{q,i,p} | \le \sum_{p=1}^n C\la _p |a_{q,i,p} | 
\le C\sum_{p=1}^n \la _p (\phi _q (a_{q,i,p} ) +1) \\
&=C(X_q (a_{q,i} ) +1) =Cu_\nu +C.
\end{align*}
Thus we have either $h(u_\nu )=0$ (if $\nu =s+1$) or $h(u_\nu ) =\ep _0 f(a)$ for
some $a\in Q_N$ such that $|a| \le Cu_\nu +C \le C|x| +9C$. Similarly, we have
either $h(u_{\nu -1} )= 0$ or $h(u_{\nu -1} )= \ep _0 f(b)$,
where $|b|\le C|x| +9C$. Therefore, $|f(a)|, |f(b)|\le \V f\V _{T_x}$ and
$$|h(x)|\le  \max (|h(u_{\nu -1} )|, |h(u_\nu )|) \le
\ep _0 \V f\V _{T_x} ,$$
proving \eqref{e7}. If $|x|\le D_t$, then $C|x|+9C\le D_{t+1}$, and thus
\eqref{e2} follows from \eqref{e7}.

\noi
{\it Proof of} \eqref{e3}: Let $x= (x_1 \stb x_n )\in Q_t$ be given, where
$t\le N$. Then $x \in \bigcup_{i\in \zn} B_q (i)$ holds for at least $m-n$ of
the indices $q$. We may assume that this is true for every $q=1\stb m-n$.

If $1\le q\le m-n$, then $x\in B_q (i)$ for a suitable $i\in \zn$. Then,
by the definition of the function $h$, we have
$h(X_q (x))=h(u_{q,i})=\ep _0  f(a_{q,i})$, and thus
$$|h(X_q (x)) -\ep _0  f(x)|=\ep _0  |f(a_{q,i})-f(x)| < \ep _0  \eta /m ,$$
since $a_{q,i} , x\in Q_N$ and $|a_{q,i} - x|<mn\de$.
Then, by $\ep _0  =1/(m-n)$ we obtain
\begin{equation}\label{e5}
  \left| f(x)-\sum_{q=1}^{m-n} h(X_q (x))\right| \le (m-n) \ep _0  \eta /m <
\ep _0  \eta  < \eta .
\end{equation}
For every $q$ we have
\begin{equation}\label{e8}
|X_q (x)|=\left| \sum_{p=1}^n \la _p  \phi _q (x_p )\right| \le
\sum_{p=1}^n  \la _p  (|x_p |+1 ) \le D_t +1 ,
\end{equation}
since $x\in Q_t$. Then by \eqref{e7} we get
$|h(X_q (x))| \le \ep _0   \V f\V _T$, where
$$T =[-C|X_q (x)|-9C, C|X_q (x)|+9C]^n .$$
Now $C|X_q (x)|+9C\le CD_t +10C=D_{t+1}$ by \eqref{e8}, and
thus $|h(X_q (x))|\le \ep _0  \V f\V _{Q_{t+1}}$. Therefore,
$$\left| \sum_{q=m-n+1}^{m} h(X_q (x))\right| \le n\ep _0  \V f\V _{Q_{t+1}} =
\ep _1 \V f\V _{Q_{t+1}} .$$
Comparing with \eqref{e5}, we obtain \eqref{e3}. 
This completes the proof of the lemma. \hfill $\square$ 

\begin{lemma}\label{l2}
There exist functions $\phi _1 \stb \phi _m \in \Phi$ such that
$(\phi _1 \stb \phi _m )\in \Om _{N,\eta} (f)$ for every $f\in C(\rn )$,
$N\in \nn$ and $\eta >0$.
\end{lemma}
\proof Let $\ff$ be a countable dense subset of $C(\rn )$, and let $W$ be the
intersection of the sets $\Om _{N,\eta} (f)$, where $f\in \ff$, $N\in \nn$
and $\eta >0$ is rational. Then $W$ is comeager in $\Phi ^m$. 
Let $(\phi_1 \stb \phi _m )\in W$ be arbitrary; we show that
$\phi_1 \stb \phi _m$ satisfy the requirements of the lemma.

Let $f\in C(\rn )$, $N\in \nn$ and $\eta >0$ be given. We may assume that
$\eta$ is rational. Since $\ff$ is dense in $C(\rn )$,
there is an $f_0 \in \ff$ such that $\V f-f_0 \V _{Q_{N+2}} <\eta /3$.
Now $(\phi _1 \stb \phi _m )\in \Om _{N,\eta /3} (f_0 )$ implies that
there exists a $h\in C(\rr )$ such that \eqref{e1}-\eqref{e3} hold true
with $f_0$ in place of $f$ and $\eta /3$ in place of $\eta$. It is clear that
$h$ satisfies \eqref{e1}. We show that it also satisfies \eqref{e2} and
\eqref{e3}.

We have
$$\V f_0 \V _{Q_t} \le\V f \V _{Q_t} + \V f-f_0 \V _{Q_t} <
\V f \V _{Q_t} +\eta /3$$
for every $t\le N+2$. This implies  
$$\V h\V _{[-D_t ,D_t ]} < \ep _0 \V f_0 \V _{Q_{t+1}} +\eta /3
< \ep _0 \V f \V _{Q_{t+1}} + \ep _0 \eta /3 +\eta /3 <
\ep _0 \V f \V _{Q_{t+1}}  +\eta$$
for every $t\le N$. If $t\ge N+1$, then taking \eqref{e1} into consideration, we
obtain
\begin{align*}
  \V h\V _{[-D_t ,D_t ]} &= \V h\V _{[-D_{N+1} ,D_{N+1} ]} < \ep _0 \V f_0 \V _{Q_{N+2}} +
  \eta /3 \\
  &\le \ep _0 \V f \V _{Q_{N+2}} +\ep _0 \eta /3 +\eta /3  \le
  \ep _0 \V f \V _{Q_{t+1}} +\eta ,
\end{align*}
proving that $h$ satisfies \eqref{e2}. As for \eqref{e3}, we have
\begin{align*}
\bigg\V f(x_1 \stb x_n ) -\sum_{q=1}^{m} h &\left( \sum_{p=1}^n \la _p  
\phi _q (x_p ) \right) \bigg\V _{Q_t} < \ep _1 \V f_0 \V _{Q_{t+1}} +\eta /3 +
\V f-f_0 \V _{Q_t} \\
&< (\ep _1 \V f \V _{Q_{t+1}} +\ep _1 \eta /3)+\eta /3 +\eta /3 \\
&<\ep _1 \V f \V _{Q_{t+1}} +\eta
\end{align*}
for every $t\le N$. \hfill $\square$

\section{Proof of Theorem \ref{t1}}
We fix a sequence $(\phi _1 \stb \phi _m )$ such that $(\phi _1 \stb \phi _m )
\in \Om _{N,\eta} (f)$ for every $f \in C(\rn )$,
$N\in \nn$ and $\eta >0$. We prove that for every bounded $f \in C(\rn )$
there exists a $g\in C(\rr )$ such that \eqref{e10} holds. We
put $X_q (x)= \sum_{p=1}^n \la _p \phi _q (x_p )$ for every $q=1\stb m$. 
\begin{lemma}\label{l4}
We have $\ep _1  (1+m\ep _0 )<1$. 
\end{lemma}
\proof
Put $x=m/n$. Then $\ep _0 =1/(m-n)=(n(x-1))\am$, $\ep _1 =1/(x-1)$,
and
$$\ep _1  (1+m\ep _0 )= (x-1)^{-1} \left( 1+\frac{x}{x-1}\right) =
(x-1)^{-2} (2x-1) <1$$
if $x> 2+\sqrt 2$. In our case $x=m/n>2+\sqrt 2$ by $m>(2+\sqrt 2 )n$,
so the inequality does hold. \hfill $\square$

We fix a number $\ep >\ep _0$ such that $\ep _1  (1+m\ep )<1$, and fix a 
positive number $\al <1/2$ such that
\begin{equation}\label{e6}
\ep _1^{1-2\al}  (1+m\ep )<1 .
\end{equation}

Let the bounded function $f_0 \in C(\rn )$ be given; we may assume that
$\V f_0 \V _{\rn} \le 1$.
We construct the sequences of functions $f_k$ and $h_k$ as follows. Let
$k\ge 0$, and suppose that $f_k \in C(\rn )$ has been defined. If $f_k =0$,
then we define $h_k =0$ and $f_{k+1} =0$. Otherwise we put
$$\eta _k = \min (\ep _1^k , (\ep -\ep _0 )\V f_k \V _{\rn} ).$$ 
Since $(\phi _1 \stb \phi _m )\in \Om _{k ,\eta _k} (f_k )$, there exists a
$h_k \in C(\rr )$ satisfying \eqref{e1}-\eqref{e3} with $h_k$ in place
of $h$, $f_k$ in place of $f$, $k$ in place of $N$ and $\eta _k$ in
place of $\eta$. Then we put $f_{k+1} =f_{k}-\sum_{q=1}^m h_{k} \circ X_q$.
In this way we define $f_k$ and $h_k$ for every $k=0,1,\ldots$.
Then we have 
\begin{equation}\label{e9}
f_0 =f_k +\sum_{q=1}^m (h_0 +\ldots +h_{k-1})\circ (X_q )
\end{equation}
for every $k$. We will show that $f_k \to 0$ 
of $\rn$, and that $\sum_{k=0}^\infty h_k$ is convergent uniformly
on every bounded subset of $\rr$.  Then \eqref{e9} will give \eqref{e10}
with $f_0$ is place of $f$, where $g=\sum_{k=0}^\infty h_k$ is continuous.

\begin{lemma}\label{l3}
We have $\V f_k \V _{\rn} \le (1+m\ep )^k$ for every $k\ge 0$.
\end{lemma}
\proof
Since $\V f_0 \V _{\rn} \le 1$ by assumption, the statement
is true for $k=0$. Let $k\ge 0$, and suppose the statement is true for $k$.
If $f_k =0$ then $f_{k+1} =0$, and the statement is true for $k+1$.
Otherwise  we have, by \eqref{e1} and \eqref{e2},
\begin{align*}
\V h_k \V _\rr & =\V h_k \V _{[-D_{k+1},D_{k+1}]} < \ep _0 \V f_k \V _{Q_{k+2}}+\eta _k\\
&\le \ep _0 \V f_k \V _{\rn} +(\ep -\ep _0 )\V f_k \V _{\rn} =
\ep  \V f_k \V _{\rn} .
\end{align*}
Thus
\begin{align*}
  |f_{k+1} (x)|&\le |f_k (x)|+\left| \sum_{q=1}^m h_k (X_q (x))\right|  
  \le \V f_k \V _{\rn}  + m\cd \ep \V f_k \V _{\rn} \\
  &=(1+m\ep )\cd \V f_k \V _{\rn} \le (1+m\ep )^{k+1}
\end{align*}
for every $x\in \rn$, proving the lemma by induction. \hfill $\square$

We put $M_{k,i}=\V f_k \V _{Q_i}$ for every $k\ge 0$ and $i\ge 1$.
\begin{lemma}\label{l5}
We have $M_{k,i} \le (k+1) \ep _1 ^{\al k-(1-\al )i}$ for every $k\ge 0$ and
$i\ge 1$.
\end{lemma}
\proof We have $\V f_0 \V _{\rn} \le 1$ by assumption. Since $\ep _1 <1$,
the statement is true for $k=0$ and for every $i\ge 1$. 

Let $k\ge 0$, and suppose the inequality holds for every $i\ge 1$.
If $1\le i\le k$, then \eqref{e3} gives
\begin{align*}
M_{k+1, i} &\le \ep _1 \V f_k \V _{Q_{i+1}}  + \ep _1 ^k 
\le \ep _1 \cd (k+1) \ep _1 ^{\al k-(1-\al )(i+1)}  + \ep _1 ^k\\
&\le (k+2) \ep _1 ^{\al (k+1)-(1-\al )i} ,
\end{align*}
since $1+\al k-(1-\al )(i+1)=\al (k+1)-(1-\al )i$ and $k\ge \al (k+1)-
(1-\al )i$.

If $i\ge k+1$ then we have, by Lemma \ref{l3} and by \eqref{e6},
$$M_{k+1, i} \le \V f_{k+1} \V _{\rn} \le (1+m\ep )^{k+1} < \ep _1 ^{(2\al -1)(k+1)}
\le  \ep _1 ^{\al (k+1)-(1-\al )i}.$$
This proves the lemma by induction. \hfill $\square$

If $i\ge 1$ is fixed, then $\V f_k \V _{Q_i} \le (k+1) \ep _1 ^{\al k-(1-\al )i}$
for every $k$. Since $0<\ep _1 <1$, we have $(k+1) \ep _1 ^{\al k} \to 0$
as $k\to \infty$. This proves that $f_k \to 0$  on $\rn$. As for the series
$\sum_{k=0}^\infty h_k$, note that for every $i\ge 1$ we have, by \eqref{e2},
$$\V h_k \V _{[-D_i ,D_i ]} \le \ep _0 \V f_k \V _{Q_{i+1}} +\ep _1^k \le M_{k,i+1}
+\ep _1 ^k \le (k+1) \ep _1 ^{\al k-(1-\al )(i+1)} +\ep_1^k$$
for every $k$. Since $\sum_{k=0}^\infty (k+1)\ep _1 ^{\al k} <\infty$
and $\sum_{k=0}^\infty \ep _1 ^k <\infty$, it follows
that the series $\sum_{k=0}^\infty h_k$ is convergent uniformly
on every bounded subset of $\rr$. As we saw above, this proves
\eqref{e10}. \hfill $\square$

\section{Sketch of proof of Theorem \ref{t2}}
Suppose $m>(2+\sqrt 2 )(2n-1)$.
We fix the vectors $\la_q =(\la _{q,1} \stb \la _{q,n} )$ $(q=1\stb m)$ with the
following properties: (i) $\la _{q,p} >0$ for every $q$ and $p$,
(ii) $\sum_{p=1}^n \la _{q,p} =1$ for every $q=1\stb m$, and (iii) any $n$-element
subset of $\{ \la _1 \stb \la _m \}$ consists of linearly independent vectors.

Note that if we select $\{ \la _1 \stb \la _{m}\}$ randomly from the simplex
$\De =\{ (x_1 \stb x_n )\colon x_1 \stb x_n  >0, \ \sum_{p=1}^n x_p =1\}$,
then almost every choice will give a system with the properties above.

It is easy to prove that there is a positive constant $C$ only depending on
the numbers $\la _{q,p}$ such that for every $x\in \rn$ we have
$|x| \le C|\lan \la _q , x\ran |$ for at least $m-n+1$ of the indices
$q=1\stb m$. We put $D=C+4$ and $Q_t =[-D^t ,D^t ]^n$ for every $t\in \nn$.
We denote $\ep _0 =(m-n+1)\am$ and $\ep _1 =(n-1)\ep _0$.

Let $\Phi$ denote the set of increasing continuous functions $f\colon \rr \to
\rr$ such that $|f(x)-x|<1$ for every $x\in \rr$. Then $\Phi$ is a $G_\de$
subset of $C(\rr )$. Therefore, $\Phi$ equipped with the subspace topology is a
Polish space. Then $\Phi ^m$, as a product space, is also Polish.

If $f\in C(\rn )$, $N\in \nn$ and $\eta >0$, then let $\Om _{N, \eta} (f)$
denote the set of those $m$-tuples $(\phi _1 \stb \phi _{m})\in \Phi ^{m}$
for which there exists a function $h\in C(\rr )$ such that
\begin{enumerate}[{\rm(i)}]
\item $h(x)=0$ if $|x|\ge D^N +2$,
\item $\V h\V _{[-D^t ,D^t ]} \le \ep _0\V f\V _{Q_{t+1}} +\eta$ for every $t\in \nn$,
and
\item $\left\V f(x_1 \stb x_n ) -\sum_{q=1}^{m} h\left( \sum_{p=1}^n \la _{q,p}
\phi _q (x_p ) \right) \right\V _{Q_t} < \ep _1 \V f\V _{Q_{t+1}} +\eta$
for every $1\le t\le N$.  
\end{enumerate}
Then $\Om _{N,\eta}(f)$ is a dense open subset of $\Phi ^m$. To prove that
$\Om _{N,\eta}(f)$ is dense we define the intervals $I_q (j)$, $K_q (i)$, the
boxes $B_{q,i}$ and the elements $a_{q,i} \in B_{q,i}$ the same way as in the proof
of Lemma \ref{l1}. Then we choose  $\phi _1 \stb \phi_m \in \Phi$
close to the given functions $\psi _q$ and such that the values
$u_{q,i} =X_q (B_q (i))$ are distinct, where $X_q =\sum_{p=1}^n \la _{q,p} \phi _q (x_p )$.

Then we define the function $h$ as follows. We put $h(u_{q,i})=
\ep _0 f(a_{q,i} )$ if $|a_{q,i}|\le C (|u_{q,i}|+2)$, and $h(u_{q,i})=0$
otherwise. We extend $h$ from the set of numbers $u_{q,i}$ to $\rr$ 
as in the proof of Lemma \ref{l1}. In the proof of (iii) we show that
the sum $\sum_{q=1}^{m} h\left( \sum_{p=1}^n \la _{q,p} \phi _q (x_p ) \right)$
contains at least $m-n+1$ terms in which $h(u_{q,i})=
\ep _0 f(a_{q,i} )$. It is this point where we use the special properties of the
vectors $\la _1 \stb \la _m$. Then the rest of the proof is the same as that
of Lemma \ref{l1}.

For the estimates of the functions $f_k$ we need $\ep _1 (1+m\ep _0 )<1$ as
in Lemma \ref{l4}, and this is why we need the bound $m>(2+\sqrt 2 )(2n-1)$.
Otherwise the conclusion of the proof is the same as in the case of
Theorem \ref{t1}.

\end{document}